# Algorithm to generate ideals in a Lie algebra of matrices at any particular characteristic with *Mathematica*

**Pablo Alberca Bjerregaard**

*Department of Applied Mathematics*
*Escuela Técnica Superior de Ingeniería Industrial*
*Campus de El Ejido*
*Universidad de Málaga*
*Málaga –Spain*
*E−Mail:pgalberca@uma.es*

**Cándido Martín González**

*Department of Algebra, Geometry and Topology*
*Facultad de Ciencias*
*Universidad de Málaga*
*Campus de Teatinos*
*Málaga –Spain*
*E−Mail:candido@apncs.cie.uma.es*

*2001.*

**Abstract.** *We present in this paper a routine which construct the ideal generated by a list of elements in a matrix Lie algebra at any particular characteristic. We have used this algorithm to analyze the problem of the simplicity of some Lie algebras (see the references). Notice that it is not included the particular construction of a basis (chosen a Lie algebra) denoted by $\{x_i\}_{i=1}^n$ .*

## ■ Previous definitions

As we have comment in the abstract, there are some commands that must be pre−definedin order to use the algorithms, listed below:

1. A basis of the Lie algebra, denoted by $x_i$ , $i = 1, \ldots, n$.
2. The dimension of the Lie algebra, denoted by $n$.
3. A function which recognize a element of the Lie algebra, denoted by **Recog**, i.e. write an element as a non−evaluatedlinear combination of the elements of the basis.

As an example, we present the Lie algebra $M_2(K)$ with char$(K) = k$. Then we define

```
x₁ = (1 0)    x₂ = (0 1)    x₃ = (0 0)    x₄ = (0 0)    n = 4;
     (0 0);        (0 0);        (1 0);        (0 1);
Recog[m_] := Module[{},
    λ₁ = m[[1, 1]]; λ₂ = m[[1, 2]]; λ₃ = m[[2, 1]];
    λ₄ = m[[2, 2]]; Sum[λᵢ * (HoldForm[x])ᵢ, {i, n}]]
```

For instance, we make

```
Recog[(a b)
      (c d)]
```

$a\, x_1 + b\, x_2 + c\, x_3 + d\, x_4$

It is also important to remark that we have not included any filter in order to check if the data is well created.

## ■ Routines

We define the Lie bracket

```
c[x_, y_] := x.y - y.x
```

and we can, for example, construct the multiplication table of the Lie algebra

```
Table[Recog[c[xᵢ, xⱼ]], {i, n}, {j, n}]
```

$$\begin{pmatrix} 0 & x_2 & -x_3 & 0 \\ -x_2 & 0 & x_1 - x_4 & x_2 \\ x_3 & x_4 - x_1 & 0 & -x_3 \\ 0 & -x_2 & x_3 & 0 \end{pmatrix}$$

The function

```
NonZero[lis_] := Module[{new = {}},
    Do[
      If[Not[lis[[i]] === 0 * lis[[i]]], AppendTo[new, lis[[i]]]],
        {i, Length[lis]}]; new]
```

works as

```
NonZero[{1, 2, 0, 3, 0, 0, 4, 5, 0}]
```

{1, 2, 3, 4, 5}

The following comands

```
Coord[y_] := Module[{coo, aa}, coo = Array[aa, {n}];
       Do[aa[i] = Coefficient[y, (HoldForm[x])ᵢ], {i, n}];
       coo]

Comb[lis_] := Module[{},
       Sum[lis[[i]] * (HoldForm[x])ᵢ, {i, n}]]
```

can be explained with the following two examples

```
Coord[Recog[( a  b
              c  d )]]
```

{*a*, *b*, *c*, *d*}

```
Comb[{1, 2, a, 4}]
```

$x_1 + 2 x_2 + a x_3 + 4 x_4$

We have then the principal algorithm. It needs two entries: a list of element (perhaps only one) and the characteristic of the base field. As a result, we have the ideal generated by the list of elements, a basis, and a list of partial results in terms of the depth of the Lie brackets.

```
IdealGenerated[lis_, char_] := Module[{k, ideal},
      ideal = Map[Comb, NonZero[
    RowReduce[Map[Coord, Map[Recog, lis]], Modulus -> char]]];
Print["Depth = ", 0, " -> ", ideal];
      dim[0] = 0; dim[1] = 1; k = 1;
      While[
   dim[k] < n && dim[k] ≠ dim[k - 1], conj = ideal; k = k + 1;
   Do[AppendTo[ideal, Recog[c[ReleaseHold[conj[[j]]], x_i]]],
    {i, n}, {j, Length[conj]}];
          ideal = Map[Comb, NonZero[
     RowReduce[Map[Coord, ideal], Modulus -> char]]];
          Print["Depth = ", k - 1, " -> ", ideal];
   dim[k] = Length[ideal]];
       Print["Ideal <", Map[Recog, lis],
   "> = ", ideal, " with ", " dimension = ",
   dim[k], " and char(K)=", char]]
```

Using our example, we can compute

```
IdealGenerated[{x_2}, 2]
```

Depth = 0 –> {$x_2$}

Depth = 1 –> {$x_1 + x_4$, $x_2$}

Depth = 2 –> {$x_1 + x_4$, $x_2$}

Ideal <{$x_2$}> = {$x_1 + x_4$, $x_2$} with dimension = 2 and char(K)=2

```
IdealGenerated[{x_1}, 2]
```

Depth = 0 –> {$x_1$}

Depth = 1 –> {$x_1$, $x_2$, $x_3$}

Depth = 2 –> {$x_1$, $x_2$, $x_3$, $x_4$}

Ideal <{$x_1$}> = {$x_1$, $x_2$, $x_3$, $x_4$} with dimension = 4 and char(K)=2

```
IdealGenerated[{x_1, x_2}, 2]
```

Depth = 0 –> {$x_1$, $x_2$}

Depth = 1 –> {$x_1$, $x_2$, $x_3$, $x_4$}

Ideal <{$x_1$, $x_2$}> = {$x_1$, $x_2$, $x_3$, $x_4$} with dimension = 4 and char(K)=2

**IdealGenerated[{x$_3$, x$_3$ - x$_1$}, 3]**

Depth = 0 –> {$x_1$, $x_3$}

Depth = 1 –> {$x_1$, $x_2$, $x_3$, $x_4$}

Ideal <{$x_3$, $x_3 - x_1$}> = {$x_1$, $x_2$, $x_3$, $x_4$} with dimension = 4 and char(K)=3

**IdealGenerated[{x$_2$}, 3]**

Depth = 0 –> {$x_2$}

Depth = 1 –> {$x_1 + 2 x_4$, $x_2$}

Depth = 2 –> {$x_1 + 2 x_4$, $x_2$, $x_3$}

Depth = 3 –> {$x_1 + 2 x_4$, $x_2$, $x_3$}

Ideal <{$x_2$}> = {$x_1 + 2 x_4$, $x_2$, $x_3$} with dimension = 3 and char(K)=3

# References

Alberca, P., Martín González, C., Elduque Palomo, A. and Navarro Márquez, F.J. *On the Cartan−Jacobson Theorem*. Journal of Algebra. In press. 2002.

Alberca, P. and Martín González, C. *On derivation algebras of some exceptional Jordan algebras in characteristic two*. Preprint. 2001.

Alberca, P. *On the Cartan−Jacobson and Chevalley−Schafer theorems. Ph. Thesis.* 2001.

Jacobson, N. *Lie algebras*. Dover Publications, Inc. 1979.

Kaufman, S. *Mathematica as a tool. An introduction with practical examples*. Birkäuser. 1994.

Postnikov, M. *Lecons de géometrie: groupes et algèbres de Lie*. Mir. 1990.

Wolfram, S. *Mathematica. A system for doing mathematics by computer*. Addison−Wesley. 1988.